\documentclass[11pt]{article}
\begin{document}
\topmargin0.5cm \headheight0.5cm \headsep0.5cm \topskip0.5cm
\textheight21cm \textwidth14.0cm \footskip2cm \oddsidemargin1.3cm
\evensidemargin1.3cm


\newtheorem{thm}{Theorem}[section]


\newtheorem{lem}[thm]{Lemma}
\newtheorem{cor}[thm]{Corollary}
\newtheorem{ex}[thm]{Example}
\newtheorem{prop}[thm]{Proposition}
\newtheorem{remark}[thm]{Remark}
\newtheorem{coun}[thm]{Counterexample}
\newtheorem{defn}[thm]{Definition}

\newcommand\ack{\section*{Acknowledgement.}}

\newcommand{\etal}{{\it et al. }}

\newcommand{\Fe}{F_{\epsilon}}
\newcommand{\Fee}{F_{\epsilon^{(1)}}}
\newcommand{\fe}{f_{\epsilon}}

\newcommand{\bbP}{{\rm I\hspace{-0.8mm}P}}
\newcommand{\bbE}{{\rm I\hspace{-0.8mm}E}}
\newcommand{\bbF}{{\rm I\hspace{-0.8mm}F}}
\newcommand{\bbI}{{\rm I\hspace{-0.8mm}I}}
\newcommand{\bbR}{{\rm I\hspace{-0.8mm}R}}
\newcommand{\bbRp}{{\rm I\hspace{-0.8mm}R}_+}
\newcommand{\bbN}{{\rm I\hspace{-0.8mm}N}}
\newcommand{\bbC}{{\rm C\hspace{-2.2mm}|\hspace{1.2mm}}}
\newcommand{\bbD}{{\rm I\hspace{-0.8mm}D}}
\newcommand{\bbQ}{\bf Q}
\newcommand{\bbZ}{{\rm \rlap Z\kern 2.2pt Z}}
\newcommand{\bbK}{{\rm I\hspace{-0.8mm}K}}

\newcommand{\matP}{{\bbP}}
\newcommand{\mattildeP}{\tilde{\bbP}}
\newcommand{\matPN}[1]{{\bbP}_{#1}}
\newcommand{\matPP}[1]{{\bbP}_{#1}^0}
\newcommand{\matE}{{\bbE}}
\newcommand{\mattildeE}{\tilde{\bbE}}
\newcommand{\matEP}[1]{{\bbE}_{#1}^0}
\newcommand{\matF}{{\bbF}}
\newcommand{\matR}{{\bbR}}
\newcommand{\matRp}{{\bbRp}}
\newcommand{\matN}{{\bbN}}
\newcommand{\matZ}{{\bbZ}}
\newcommand{\matI}{{\bbI}}
\newcommand{\matK}{{\bbK}}
\newcommand{\matQ}{{\bbQ}}
\newcommand{\matC}{{\bbC}}
\newcommand{\matD}{{\bbD}}

\newcommand{\calL}{{\cal L}}
\newcommand{\calM}{{\cal M}}
\newcommand{\calN}{{\cal N}}
\newcommand{\calF}{{\cal F}}
\newcommand{\calG}{{\cal G}}
\newcommand{\calD}{{\cal D}}
\newcommand{\calB}{{\cal B}}
\newcommand{\calH}{{\cal H}}
\newcommand{\calI}{{\cal I}}
\newcommand{\calP}{{\cal P}}
\newcommand{\calQ}{{\cal Q}}
\newcommand{\calS}{{\cal S}}
\newcommand{\calT}{{\cal T}}
\newcommand{\calC}{{\cal C}}
\newcommand{\calK}{{\cal K}}
\newcommand{\calX}{{\cal X}}
\newcommand{\cals}{{\cal S}}
\newcommand{\calE}{{\cal E}}

\newcommand{\koniecmat}{\,}

\newcommand{\eqd}{\ =_{\rm d}\ }
\newcommand{\toto}{\leftrightarrow}
\newcommand{\eqdistr}{\stackrel{\rm d}{=}}
\newcommand{\convdistr}{\stackrel{\rm d}{\rightarrow}}
\newcommand{\convweak}{{\Rightarrow}}
\newcommand{\convas}{\stackrel{\rm a.s.}{\rightarrow}}
\newcommand{\convfidi}{\stackrel{\rm fidi}{\rightarrow}}
\newcommand{\convprob}{\stackrel{p}{\rightarrow}}
\newcommand{\deff}{\stackrel{\rm def}{=}}
\newcommand{\bis}{{'}{'}}
\newcommand{\Cov}{{\rm Cov}}
\newcommand{\Var}{{\rm Var}}
\newcommand{\Exp}{{\rm E}}

\newcommand{\nd}{n^{\delta}}
\newcommand{\koniec}{\newline\vspace{3mm}\hfill $\odot$}


\title{Optimal rates in the Bahadur-Kiefer representation for GARCH sequences\protect}
\author{Rafa{\l} Kulik\\
University of Sydney and Wroc{\l}aw University\thanks{School of
Mathematics and Statistics, F07, University of Sydney, NSW 2006,
Australia and Mathematical Institute, University of Wroc{\l}aw, Pl.
Grunwaldzki 2/4, 50-384 Wroc{\l}aw, Poland, email:
rkuli@math.uni.wroc.pl} }

\maketitle


\begin{abstract}
In this paper we establish the Bahadur-Kiefer representation for
sample quantiles of GARCH sequences with optimal rates.
\end{abstract}
\noindent{\bf Keywords:} sample quantiles, GARCH sequences, Bahadur
representation \noindent{\bf Running title:} Sample quantiles and
GARCH



\section{Introduction}
Over the last years it was observed that (financial) data can be
modeled appropriately by the Autoregressive Conditionally
Heteroskedastic (ARCH) sequences, introduced by Engle
\cite{Engle-1982}. This model has been generalized later to GARCH
sequences, see e.g. \cite{Bollerslev-1986}. A GARCH$(p,q)$ process
is defined by
\begin{equation}\label{eq.garch}
X_k=\sigma_k\varepsilon_k ,
\end{equation}
\begin{equation}\label{eq.garch.2}
\sigma_k^2=\delta+\sum_{i=1}^p\beta_i\sigma_{k-i}^2+\sum_{j=1}^q\alpha_jX_{k-j}^2,
\end{equation}
where $\delta>0$, $\beta_i$, $1\le i\le p$ and $\alpha_j$, $1\le
j\le q$ are nonnegative constants. We assume that
$\{\varepsilon_i,-\infty\le i\le \infty\}$ are i.i.d. random
variables with distribution $H$. Under appropriate conditions on
coefficients $\beta_i$ and $\alpha_j$ these equations
(\ref{eq.garch}) and (\ref{eq.garch.2}) have a unique stationary
solution, see \cite{Bougerol-Picard-1992-a}. In particular, let
$$
\tau_n=(\beta_1+\alpha_1\varepsilon_n^2,\beta_2,\ldots,\beta_{q-1})\in
{\matR}^{p-1}
$$
and
$$
\xi_n=(\varepsilon_n^2,0,\ldots,0)\in {\matR}^{p-1}, \qquad
\underline{\alpha}=(\alpha_2,\ldots,\alpha_{q-1})\in {\matR}^{q-1}.
$$
Further, the $(p+q-1)\times (p+q-1)$ matrix $A_n$ is defined by
$$
A_n=\left[
\begin{array}{cccc}
\tau_n & \beta_p & \underline{\alpha} & \alpha_q\\
I_{p-1} & 0 & 0 & 0 \\
\xi_n & 0 & 0 & 0\\
0 & 0 & I_{q-2} & 0
\end{array}
\right],
$$
where $I_r$ is the identity matrix of size $r$. Let $||\cdot ||$ be
the matrix norm. Then, under the condition $\Exp
(\log^+||A_0||)<\infty$, we have
$$
\gamma=\lim_{n\to\infty}\frac{1}{n}\log||A_0A_{-1}\cdots A_{-n}||
$$
almost surely. Bougerol and Picard \cite{Bougerol-Picard-1992-a}
showed
that the unique stationary solution exist if and only if $\gamma<0$.\\

The aim of this paper is to obtain the Bahadur-Kiefer representation
for sample quantiles in case of stationary GARCH sequences. To state
our results, assume that the stationary sequence $\{X_i,i\ge 1\}$
has marginal distribution function $F(x)=P(X_1\le x)$ and a density
$f=F'$. Given the sample $X_1,\ldots,X_n$, let
$F_n(x)=\frac{1}{n}\sum_{i=1}^n1_{\{X_i\le x\}}$ be its
corresponding empirical distribution function. Let
$X_{1:n}\le\cdots\le X_{n:n}$ be the order statistics. The empirical
quantile function $Q_n(y)$, $y\in (0,1)$ is defined as
$Q_n(y)=X_{k:n}$ if $\frac{k-1}{n}\le p\le \frac{k}{n}$. Further,
let
$$
\beta_{n}(x)=n^{1/2}(F_n(x)-F(x)),\ \ \ x\in {\matR}\; ,
$$
$$
q_n(y)=n^{1/2}(Q(y)-Q_n(y)),\ \ \ y\in (0,1)\; ,
$$
be the empirical and the quantile processes, respectively. Since $F$
is continuous, we may define $U_i=F(X_i)$, $i\ge 1$. Let
$E_n(x)=\frac{1}{n}\sum_{i=1}^n1_{\{E_i\le x\}}=F_n(Q(x))$ and
$G_n(y)$ be the corresponding uniform empirical distribution and
uniform empirical quantile functions. Let
$$
\alpha_{n}(x)=n^{1/2}(E_n(x)-x),\ \ \ x\in [0,1]\; ,
$$
$$
\gamma_n(y)=n^{1/2}(y-G_n(y)),\ \ \ y\in (0,1)\; ,
$$
be the corresponding uniform empirical and uniform quantile
processes.\\

Assume for a while that $X_i$, $i\ge 1$ are i.i.d. Kiefer
\cite{Kiefer-1970}, continuing his previous research
(\cite{Kiefer-1967}) and that of Bahadur (\cite{Bahadur-1966})
obtained, in particular, the following {\it Bahadur-Kiefer
representation}:
\begin{equation}\label{Kiefer}
\sup_{y\in [0,1]}\left|u_n(y)-\alpha_n(y)\right|=:R_n
\end{equation}
and
\begin{equation}\label{Kiefer-rate}
R_n=O_{a.s}(n^{-1/4}(\log n)^{1/2}(\log\log n)^{1/4}).
\end{equation}
The above rate is exact and constants can be given, see e.g.
\cite{Deheuvels-Mason-1990}. As for the general quantile processes
the best available result is due to Cs\"{o}rg\H{o} and
R\'{e}v\'{e}sz, \cite{Csorgo-Revesz-1978}. Under appropriate
conditions (so called Cs\"{o}rg\H{o}-R\'{e}v\'{e}sz conditions) on
the distribution $F$, we have
\begin{equation}\label{CR}
\sup_{y\in (0,1)}|f(Q(y))q_n(y)-\alpha_n(y)|=O_{a.s}(n^{-1/4}(\log
n)^{1/2}(\log\log n)^{1/4}).
\end{equation}
We refer to
\cite{Csorgo-1983-LN} as well as to
\cite{Csorgo-Szyszkowicz-1998} for more discussion.\\

For weakly dependent random variables some results are available.
For mixing sequences the best possible (in terms of the rates for
$R_n$) results are included in \cite{Babu-Singh-1978}. In
particular, for a class of $\phi$-mixing sequences, they obtained
(\ref{Kiefer}) with the optimal rate (\ref{Kiefer-rate}). Those
results were improved in \cite{Fotopoulos-Ahn-1994} and
\cite{Yu-1996} in terms of less restrictive mixing rates and
Cs\"{o}rg\H{o}-R\'{e}v\'{e}sz conditions. However, the rate was
$R_n=O_{a.s}((\log n)^{-\lambda})$ with some $\lambda>0$. This rate
is much worse compared to the optimal one in (\ref{Kiefer-rate}).

However, mixing is rather hard (if possible at all) to verify and
requires some additional regularity assumptions. In particular, for
linear processes
\begin{equation}\label{linear-model}
X_i=\sum_{k=0}^{\infty}c_k\epsilon_{i-k}, \end{equation} in order to
obtain strong mixing both regularity assumptions on a density of
$\epsilon_1$ and some constrains on $c_k$'s are required (cf.
\cite{Doukhan-1984}). If $c_k$ decay exponentially fast and some
regularity assumptions hold, then we are able to establish the
strong mixing with geometric rates. However, even in this case we do
not attain the optimal rate in the Bahadur-Kiefer representation,
see \cite{Babu-Singh-1978}. To overcome such problems, Ho, Hsing,
Mielniczuk and Wu (see \cite{Ho-Hsing-1997},
\cite{Wu-Mielniczuk-2002}, \cite{Hsing-Wu-2004}, \cite{Wu-2005})
developed a martingale based methods, which leads to optimal or
almost optimal results, especially in a context of weak convergence.
Based on this method and restricting to an interval $[y_0,y_1]$,
$0<y_0<y_1<1$, Wu \cite{Wu-2005} obtained for a class of linear
processes as well as for a class of weakly dependent sequences
satisfying a {\it geometric moment contraction} assumption, almost
optimal rates in the Bahadur-Kiefer representation for linear
processes.

Coming back to the GARCH processes, we note that under appropriate
conditions on moments and a density $h$ of $H(x)=P(\varepsilon_1\le
x)$, the sequence $\{X_i^2\}$, and consequently $\{X_i\}$, is
strongly mixing with a geometric rate. However, in view of
\cite{Babu-Singh-1978}, the rates are not optimal. Therefore, in
this paper, we shall obtain the Bahadur-Kiefer representation for
GARCH sequences with the optimal rate. We note that to do this we
need, essentially, two types of results. First, we need an uniform
law of the iterated logarithm (ULIL) for the empirical process based
on GARCH sequence. Second, we need to control increments of the
empirical process. The ULIL will be a consequence of Berkes and
Horv\'{a}th \cite{Berkes-Horvath-2001} strong approximation result.
Increments will be controlled using the martingale approximation as
introduced in \cite{Wu-2005}. We note, that the situation is
complete different compared to linear processes, namely, in this
case, it is easy to control increments, see \cite[Section
6.3]{Wu-2005}. In case of GARCH processes, due to
the lack of linearity, the situation becomes much more involved.\\

We should also mention (personal communication with Wei Biao Wu)
that for GARCH sequences, the Bahadur representation with the
optimal rates can be obtained using the same technique as in
\cite[Theorem 1]{Wu-2005}.\\

In what follows, $C$ will denote a generic constant which may be
different at each of its appearances. Also, we write $a_n\sim b_n$
if $\lim_{n\to\infty}a_n/b_n=1$. For any stationary sequence
$\{Z_i,i\ge 1\}$ of random variables, $Z$ will be a random variable
with the same distribution as $Z_1$.\\

Throughout the paper we shall use:
$$
b_n=n^{-3/4}(\log n)^{1/2}(\log\log n)^{1/4}
$$
and
$$
\lambda_n=n^{-1/2}(2\log\log n)^{1/2}.
$$
For any function $h(x)$ defined on ${\matR}$ we write for $x<y$,
$h(x,y):=h(y)-h(x)$.
\section{Statement of results}
Note that $F(x)=P(X\le x)=\Exp H(x/\sigma)$. Thus,
$f(x)=\Exp\sigma^{-1}h(x/\sigma)$ and
$f'(x)=-\Exp\sigma^{-2}h'(x/\sigma)$. Consequently, via $\sigma\ge
\delta$, if
$$
{\rm (H1)}\quad \sup_{x\in {\matR}}|h'(x)|<\infty, \qquad {\rm (H2)}
\quad\inf_{x\in {\matR}}h(x)>0.
$$
then also
$$
{\rm (K1)}\quad \sup_{y\in (0,1)}|f'(Q(y))|<\infty, \qquad {\rm
(K2)} \quad\inf_{y\in (0,1)}f(Q(y))>0.
$$
The main result of this paper is following.
\begin{thm}\label{th-main}
Consider the stationary GARCH model. Assume that
\begin{equation}\label{eq-fourth-moment}
\Exp \varepsilon^4<\infty .
\end{equation}
Assume (H1) and (H2). Then $$\sup_{y\in
(0,1)}|f(Q(y))q_n(y)-\alpha_n(y)|=O_{a.s}(n^{-1/4}(\log
n)^{1/2}(\log\log n)^{1/4}).
$$
\end{thm}
\begin{remark}{\rm
For a general $F$, the above result holds if one restricts to
$[a,b]$, $0<a<b<1$. To have estimates on $(0,1)$ one could follow
the path of \cite{Csorgo-Revesz-1978} or like in \cite{Yu-1996}.
However, this is out of scope of our paper. Our main concern was to
obtain the optimal rates. }
\end{remark}
\section{Proof of the Theorem}
\subsection{Variance bound}
Let ${\cal F}_i=\sigma(\varepsilon_i,\varepsilon_{i-1},\ldots,)$.
Write
\begin{eqnarray*}
\lefteqn{\frac{1}{n}\sum_{i=1}^n(1_{\{X_i\le
x\}}-F(x))=}\\
&= & \frac{1}{n}\sum_{i=1}^n\left(1_{\{X_i\le x\}}-\Exp (1_{\{X_i\le
x\}}|{\cal F}_{i-1})\right)+\frac{1}{n}\sum_{i=1}^n\left(\Exp (1_{\{X_i\le x\}}|{\cal F}_{i-1})-F(x)\right)\\
&= & \frac{1}{n}\sum_{i=1}^n\left(1_{\{X_i\le x\}}-\Exp (1_{\{X_i\le
x\}}|{\cal F}_{i-1})\right)+\frac{1}{n}\sum_{i=1}^nY_i(x)=:
M_n(x)+N_n(x).
\end{eqnarray*}
Then $nM_n(x)$, $n\ge 1$, is a martingale. Also, since $\sigma_i$ is
${\cal F}_i$-measurable and $\varepsilon_i$ is independent of ${\cal
F}_i$ one has $\Exp (1_{\{X_i\le x\}}|{\cal
F}_{i-1})=H(x/\sigma_i)$.
\begin{prop}
Under the conditions of Theorem \ref{th-main},
$$
||nN_n(x,y)||_2^2\le D_0n|y-x|
$$
with some finite and positive constant $D_0$.
\end{prop}
The proof of this will be divided into several steps. First, we have
\begin{eqnarray*}
\lefteqn{\left\|nN_n(x,y)\right\|_2^2=\left\|\sum_{i=1}^n\left(H(x/\sigma_i,y/\sigma_i)-\Exp
H(x/\sigma_i,y/\sigma_i)
\right)\right\|_2^2}\\
&=&\left\|\sum_{i=1}^n\left(\int_{x/\sigma_i}^{y/\sigma_i}h(u)du-\Exp
\int_{x/\sigma_i}^{y/\sigma_i}h(u)du\right)\right\|_2^2\\
 &= &
\left\|\sum_{i=1}^n\left(\frac{1}{\sigma_i}-\Exp\frac{1}{\sigma_i}\right)\left(\int_{x}^{y}h(u)du\right)\right\|_2^2\le
C|y-x|\times
\left\|\sum_{i=1}^n\left(\frac{1}{\sigma_i}-\Exp\frac{1}{\sigma_i}\right)\right\|_2^2
\end{eqnarray*}
Thus, it suffices to show
\begin{prop}\label{prop-On}
Under conditions of Theorem \ref{th-main},
$$
\left\|\sum_{i=1}^n\left(\frac{1}{\sigma_i}-\Exp\frac{1}{\sigma_i}\right)\right\|_2^2=O(n).
$$
\end{prop}
Let us start with several results, which will be needed in the
sequel. First, under conditions of Theorem \ref{th-main},
\begin{equation}\label{cnvergence-covariances}
\sum_{i=-\infty}^{\infty}\Cov (X_0^2,X_i^2)<\infty ,
\end{equation}
(see \cite[Theorem 2.5]{Berkes-Horvath-Kokoszka-2004}). Next, the
conditional variances $\sigma_k$ can be represented in terms of
ARCH($\infty$) sequences:
\begin{equation}\label{archinfty}
\sigma_k^2=a+\sum_{i=1}^{\infty}b_iX_{k-i}^2
\end{equation}
with nonnegative and summable coefficients $b_i$, see
\cite{Kazakevicius-Leipus} and references therein. Therefore, by
(\ref{cnvergence-covariances}), (\ref{archinfty}) and stationarity
\begin{equation}\label{cov-sigma} \sum_{k=0}^{\infty}\Cov
(\sigma_0^2,\sigma_k^2)=\sum_{i=1}^{\infty}\sum_{j=1}^{\infty}b_ib_j\sum_{k=0}^{\infty}\Cov
(X_{-i}^2,X_{k-j}^2)<\infty .
\end{equation}
The proof of Proposition \ref{prop-On} will use a particular type of
dependence structure, the so called association.
\subsubsection{Association} Random vectors $R_1,\ldots,R_k$ with
values in ${\matR}^d$ are associated if for all coordinatewise
nondecreasing functions $f,g:{\matR}^{kd}\to {\matR}$ we have
\begin{equation}\label{def.assoc}
\Cov(f(R_1,\ldots,R_k),g(R_1,\ldots,R_k))\ge 0\;
\end{equation}
whenever the covariance is defined. A sequence $\{R_k\}_{k\ge 1}$ is
associated if (\ref{def.assoc}) holds for all $k\ge 1$.

This concept of weak dependence was introduced (in the scalar case)
in \cite{Esary-Proschan-Walkup-1967} and has been widely studied
since then, see e.g. \cite{Yu-1995} and references therein. We shall
need the following properties:
\begin{itemize}
\item[{\rm (A1)}] Independent random variables are associated.
\item[{\rm (A2)}] Increasing transforms of associated random variables are
associated, i.e. if $(Z_1,\ldots,Z_m)$ are associated and $f_i$,
$i=1,\ldots,k$ are coordinatewise increasing, then
$R_1=f_1(Z_1,\ldots,Z_m),\ldots, R_k=f_k(Z_1,\ldots,Z_m)$ are
associated.
\item[{\rm (A3)}] A subset of associated random variables is
associated.
\item[{\rm (A4)}] If $T_1^{(k)},\ldots,T_n^{(k)}$ are associated for
each $k$ and $(T_1^{(k)},\ldots,T_n^{(k)})\to (T_1,\ldots,T_n)$ in
distribution as $k\to \infty$, then $T_1,\ldots,T_n$ are associated.
\end{itemize}
\begin{lem}\label{lem.association}
Under the conditions of Theorem \ref{th-main}, the sequence
$\{X_k^2,-\infty<k<\infty\}$ is associated.
\end{lem}
{\it Proof.} The solution to (\ref{eq.garch}) is given by $ {\bf
Y}_n=B+\sum_{k=0}^{\infty}A_nA_{n-1}\cdots A_{n-k}B, $ where $${\bf
Y}_n:=(Y_n(1),\ldots,Y_n(p+q-1))=(\sigma_{n-1}^2,\ldots,\sigma_{n-p+2}^2,X_n^2,\ldots,X_{n-q+2}^2).$$
Then $X_n^2$ is identical with the $(p+1)$th coordinate of ${\bf
Y}_n$. Define ${\bf Y}_{n,{(m)}}=B+\sum_{k=0}^{m}A_nA_{n-1}\cdots
A_{n-k}B$ and $X_{n,{(m)}}^2$ respectively. Since the expression on
the right hand side in the preceding summation involves the
increasing transformation of i.i.d. random variables, then by (A1)
and (A2), $Y_{n,{(m)}}(1),\ldots,Y_{n,{(m)}}(p+q-1)$, the
coordinates of ${\bf Y}_{n,{(m)}}$ are associated for each fixed $m$
and $n$. Further, for arbitrary but fixed $k$, a family
$Y_{n,{(m)}}(1),\ldots,Y_{n,{(m)}}(p+q-1)$, $n=1,\ldots,k$, is
associated for each $m$. Therefore, by (A3), for arbitrary $k$ the
random variables $X_{1,{(m)}},\ldots,X_{{k},{(m)}}$ are associated.
Now, since for each $i$,  $X_{i,{(m)}}\to X_i$ almost surely as
$m\to\infty$ (see e.g. \cite[Lemma 3.2]{Berkes-Horvath-2001}), the
result follows by (A4).\koniec

In view of (\ref{archinfty}), Lemma \ref{lem.association} and using
the same truncation argument we obtain the following corollary.
\begin{cor}\label{cor.association.sigma}
Under the conditions of Theorem \ref{th-main}, the sequence
$\{\sigma_k^2,-\infty<k<\infty\}$ is associated.
\end{cor}
\subsubsection{Proof of Proposition \ref{prop-On}} Let
$g(x)=x^{-1/2}$. Let $||g||_{[\delta,\infty)}=\sup_{x\in
[\delta,\infty)}|g(x)|$. Then, for any $(\delta>0)$,
$||g'||_{[\delta,\infty)}\le (2\delta^{3/2})^{-1}$. Therefore, we
have
\begin{eqnarray*}
\lefteqn{\left\|\sum_{i=1}^n\left(\frac{1}{\sigma_i}-\Exp\frac{1}{\sigma_i}\right)\right\|_2^2=\left\|\sum_{i=1}^n\left(g(\sigma_i^2)-\Exp
g(\sigma_i^2)\right)\right\|_2^2 }\\
&=& \sum_{i=1}^n\Var g(\sigma_i^2)+2\sum_{i<j}\Cov
(g(\sigma_i^2),g(\sigma_j^2))\le n\Var
g(\sigma_1^2)+2||g'||_{[\delta,\infty)}^2\sum_{i<j}\Cov
(\sigma_i^2,\sigma_j^2),
\end{eqnarray*}
where the last estimate follows by Corollary
\ref{cor.association.sigma} and using standard moment bounds for
associated sequences, see \cite[Lemma 3.1]{Birkel-1988}. Now, the
result follows by (\ref{cov-sigma}) and the Kronecker lemma. \koniec
\subsection{Exponential inequality}
From Proposition \ref{prop-On} and the martingale approximation,
\begin{equation}\label{eq-variance-bound}
\left\|\sum_{i=1}^nX_i(x,y)\right\|_2^2\le D_1n(y-x)
\end{equation}
with some $D_1>0$.

Choose an arbitrary $\rho\in (0,\frac{1}{4})$. As in the proof of
Lemma \ref{lem.association}, consider ${\bf Y}_i$ and ${\bf \hat
Y}_i:={\bf Y}_{i,([i^{\rho}])}$, $[\cdot]$ being the integer part.
Define $\hat X_i$ in such the way that $\hat X_i^2$ agrees with
$(p+1)$th coordinate of ${\bf \hat Y}_i$ and the sign is that of
$\varepsilon_i$. Let
$$
X_i(x)=1_{\{X_i\le x\}}-F(x),\qquad \hat X_i(x)=1_{\{\hat X_i\le
x\}}-\Exp 1_{\{\hat X_i\le x\}}.
$$
Choose an arbitrary $\nu>8$. Choose $\mu$ so big so that
$\rho(\mu-2)/(4\nu)>1$. By (\ref{eq-fourth-moment}) we have $\Exp
(\log ^+|\varepsilon|)^{\mu}<\infty$. Consequently, via Remark 1.2
in \cite{Berkes-Horvath-2001}, their condition (1.6) is fulfilled.
Then, as in \cite[Lemma 2.4]{Berkes-Horvath-2001},
$$
P(|X_i-\hat X_i|>Ci^{-\rho(\mu-2)/4})\le Ci^{-\rho(\mu-2)/4}
$$
and consequently (cf. \cite[Lemma 2.5]{Berkes-Horvath-2001} with
$\theta=1$ by differentiability of $H$)
$$
P(X_i(x,y)\not= \hat X_i(x,y))\le Ci^{-\rho(\mu-2)/4} .
$$
Therefore, with $\nu$ as above
\begin{eqnarray}
\lefteqn{\hspace*{-3cm}P\left(\sum_{i=1}^n(X_i(x,y)-\hat
X_i(x,y))>z\right)\le \frac{\left(\sum_{i=1}^n\left(\Exp
|X_i(x,y)-\hat
X_i(x,y)|^{\nu}\right)^{1/\nu}\right)^{\nu}}{z^{\nu}}}\nonumber\\
&\le &
\frac{\left(\sum_{i=1}^ni^{-\rho(\mu-2)/(4\nu)}\right)^{\nu}}{z^{\nu}}\le
Cz^{-\nu}.\label{eq-step1}
\end{eqnarray}
Further, in view of (\ref{eq-variance-bound}) and as in
(\ref{eq-step1}),
\begin{equation}\label{eq-variance-bound-truncated}
\left\|\sum_{i=1}^n\hat X_i(x,y)\right\|_2^2\le D_1n(y-x)
\end{equation}
if $n(y-x)\to\infty$. Thus, we formulate all results below under
this constrain.
\begin{lem}\label{eq.exponential.inequality}
Under the conditions of Theorem \ref{th-main}, for any $z>0$
$$
{\matP}\left(\left|\sum_{i=1}^nX_i(x,y)\right|>z\right)\le
C_1z^{-\nu}+C_2\exp(-C_3z^2/(n(y-x)))+C_4\exp(-C_5z/n^{\rho}),
$$
where $C_1,C_2,C_3,C_4,C_5$ are positive constants and $\nu>8$.
\end{lem}
{\it Proof.} From (\ref{eq-step1}) and Markov inequality one gets
$${\matP} (|nN_n(x,y)|>z)\le C_1z^{-\nu}+{\matP}(|n\hat
N_n(x,y)|>z/2).
$$
To obtain the bound for the second part, divide $[1,n]$ into blocks
$I_1,J_1,I_2,J_2,\ldots,I_M,J_M$ with the same length $n^{\rho}$,
$\rho$ as above. Thus, $M\sim n^{1-\rho}$. Let $\hat U_k=\sum_{i\in
I_k}\hat X_i(x,y)$, $\hat V_k=\sum_{i\in J_k}\hat X_i(x,y)$ and
$n\hat N_n^{(1)}=\sum_{k=1}^M\hat U_k$, $n\hat
N_n^{(2)}=\sum_{k=1}^M\hat V_k$. Both $(\hat U_1,\ldots,\hat U_M)$
and $(\hat V_1,\ldots,\hat V_M)$ are vectors of independent random
variables. Also, $\max_{k=1,\ldots,M}\hat U_k\le [I_k]\le
Cn^{\rho}$. The relation (\ref{eq-variance-bound-truncated}) yields
$||U_k||_2^2\le D_1n^{\rho}(y-x)$. Applying the result in \cite[p.
293]{Petrov-1975} to the centered sequence $\hat U_1,\ldots,\hat
U_M$ with $M_n=n^{\rho}$, $B_n=D_1n(y-x)$ one obtains
$$
{\matP}(|n\hat N_n^{(1)}|>z)\le \exp
(-Cz^2/(4D_1n(y-x)^2))+\exp(-Cz/(4n^{\rho})) .
$$
The same applies to ${\matP}(|n\hat N_n^{(2)}|>z)$ and hence the
result follows. \koniec
\begin{remark}{\rm
Note that the proof of Lemma 2.8 in \cite{Berkes-Horvath-2001} has a
gap. The result from \cite{Petrov-1975} is applied to a non-centered
sequence.}
\end{remark}

\subsection{Almost sure behavior of increments of the empirical process}
Recall that $X_i(x)=1_{\{X_i\le x\}}-F(x)$ and let
$U_i(x)=1_{\{U_i\le x\}}-x$. Then
$F_n(x)-F(x)=\frac{1}{n}\sum_{i=1}^nX_i(x)$ and
$E_n(x)-x=\frac{1}{n}\sum_{i=1}^nU_i(x)$. Let $b_n^*=n^{-1/2}b_n$.
\begin{lem}Under conditions of Theorem \ref{th-main}
one has
\begin{equation}\label{modulus-of-continuity-uniform-2}
\sup_{x,y\in {\matR},|x-y|\le
\lambda_n}|F_n(y)-F_n(x)-(F(y)-F(x))|=O_{a.s}(b_n^*).
\end{equation}
\end{lem}
Let $|y-x|<1$. Since $F$ is differentiable with the strictly
positive derivative we obtain that for a sufficiently small and
positive $h$ and $|y-x|<h$ one has $C_6(y-x)\le F(x,y)\le C_7(y-x)$,
$C_6,C_7>0$. Thus, with $C_8=C_3C_6/(4D_1)$,
\begin{eqnarray*}
\lefteqn{{\matP}\left(\left|\sum_{i=1}^nX_i(x,y)\right|>z\right)}\\
&\le &
C_1z^{-\nu}+C_2\exp(-C_8z^2/(nF(x,y)))+C_4\exp(-C_5z/n^{\rho}).
\end{eqnarray*}
Substituting $U_i=F(X_i)$, $u=Q(x)$, $v=Q(y)$ one obtains
\begin{eqnarray}
\lefteqn{{\matP}\left(\left|\sum_{i=1}^nU_i(u,v)\right|>C_8z\right)}\nonumber\\
&\le &
C_1z^{-\nu}+C_2\exp(-C_8z^2/(n(v-u)))+C_4\exp(-C_5z/n^{\rho})\label{step20}.
\end{eqnarray}
First, we show that
\begin{equation}\label{modulus-of-continuity-uniform-1}
\sup_{u\in [0,1]}\sup_{v\in[0,1],|u-v|\le
\lambda_n}|E_n(u)-E_n(v)-(u-v)|=O_{a.s}(b_n^*).
\end{equation}
Let $e_n=[1/(b_n^*)]+1$, $d_n=[\lambda_n/(b_n^*)]+1$. Note that both
$e_n$, $d_n$ diverge to $\infty$. Note that the latter expression is
bounded by
$$
\max_{|j|\le d_n}\{\max_{i\le
e_n}|E_n((i+j)b_n^*)-E_n(ib_n^*)-jb_n^*|\}+2b_n^*.
$$
Let $D>2^{1/2}(C_8)^{-1/2}$. Using (\ref{step20}) with $z=Dnb_n^*$,
bearing in mind $\rho\in (0,\frac{1}{4})$ and noting that
$e_nd_n=O(n)$,
\begin{eqnarray*}
\lefteqn{P(\max_{|j|\le d_n}\{\max_{i\le
c_n}|E_n((i+j)b_n)-E_n(ib_n)-jb_n|\}>Dnb_n^*)}\nonumber\\
&\le & e_nd_n\sup_{u\in [0,1]}\sup_{v\in [0,1],|u-v|<\lambda_n}P(|E_n(u)-E_n(v)-(u-v)|>Dnb_n^*)\nonumber\\
&=&O(n(nb_n^*)^{-\nu})+O\left(nn^{-C_8D^2}\right)+
O(n\exp(-Cn^{\omega})),
\end{eqnarray*}
with $\omega>0$. By the choice of $D$ and $\nu>8$, the bound is
summable and thus (\ref{modulus-of-continuity-uniform-1}) follows by
the Borel-Cantelli lemma.

Substituting $u=F(x)$ and $v=F(y)$ into
(\ref{modulus-of-continuity-uniform-1}) we obtain the result.
\koniec

\subsection{Conclusion of the proof of Theorem \ref{th-main}}
By Lemma \ref{eq.exponential.inequality} we have $$ \sup_{x\in
{\matR}}\sup_{|x-y|\le
\lambda_n}|\beta_n(x)-\beta_n(y)|=O_{a.s}(b_n).
$$
From \cite{Berkes-Horvath-2001} one has
\begin{equation}\label{eq-lil-empirical}
\limsup_{n\to\infty}(\log\log
n)^{-1/2}\sup_{-\infty<x<\infty}|\beta_n(x)|=C \qquad \mbox{\rm
almost surely}.
\end{equation}
If (\ref{eq-lil-empirical}) holds then taking $U_i=F(X_i)$ one
obtains the same result for the uniform empirical process. Thus, via
$\sup_{y\in (0,1)}|\gamma_n(y)|=\sup_{y\in (0,1)}|\alpha_n(y)|$,
\begin{equation}\label{eq-lil-uniform-quantile}
\limsup_{n\to\infty}(\log\log
n)^{-1/2}\sup_{-\infty<x<\infty}|\gamma_n(x)|=C \qquad \mbox{\rm
almost surely}.
\end{equation}
Further, if (H1)-(H2) are fulfilled,
\begin{equation}\label{eq-lil-quantile}
\limsup_{n\to\infty}(\log\log
n)^{-1/2}\sup_{-\infty<x<\infty}|q_n(x)|=C \qquad \mbox{\rm almost
surely}.
\end{equation}
On account of (\ref{eq-lil-quantile}) we have
$$\sup_{y\in (0,1)}|\beta_n(Q_n(y))-\beta_n(Q(y))|=O_{a.s}(b_n).
$$
 Equivalently,
$$
n^{1/2}\sup_{y\in
(0,1)}|F_n(Q_n(y))-F(Q_n(y))-(F_n(Q(y))-F(Q(y)))|=O_{a.s}(b_n).
$$
Since $|F_n(Q_n(y))-F(Q(y))|\le 1/n$ one obtains
$$
n^{1/2}\sup_{y\in
(0,1)}|F(Q_n(y))-F(Q(y))-(F(Q(y))-F_n(Q(y)))|=O_{a.s}(b_n).
$$
Set $\Delta_{n,y}=Q_n(y)-Q(y)$. Using the Taylor's expansion
$F(Q_n(y))=F(Q(y))+f(Q(y))\Delta_{n,y}+O_{a.s}(\Delta_{n,y}^2)$ we
finish the proof of Theorem \ref{th-main}. \koniec

\ack This work was done during my stay at Carleton University. I am
thankful to Professors Barbara Szyszkowicz and Miklos Cs\"{o}rg\H{o}
for the support and helpful remarks.


\begin{thebibliography}{99}
\bibitem{Babu-Singh-1978}
Babu, G.J. and Singh, K. (1978). On deviatios between empirical and
qunatile processes for mixing random variables. {\it Journal of
Multivariate Analysis} {\bf 8}, 532--549.
\bibitem{Bahadur-1966}
Bahadur, R.R. (1966). A note on quantiles in large samples. {\it
Ann. Math. Statist.} {\bf 37}, 577--580.
\bibitem{Bollerslev-1986}
Bollerslev, T. (1986). Generalized autoregressive conditional
heteroskedasticity. {\it J. Econometrics} {\bf 31}, 307--327.
\bibitem{Berkes-Horvath-2001}
Berkes, I. and Horv\'{a}th, L. (2001). Strong approximation of the
empirical process of GARCH sequences. {\it Ann. Appl. Probab.} {\bf
11}, 789--809.
\bibitem{Berkes-Horvath-Kokoszka-2004}
Berkes, I., Horv\'{a}th, L. and Kokoszka, P. (2004). Probabilistic
and statistical properties of GARCH processes. {\it Asymptotic
methods in stochastics}, 409--429, Fields Inst. Commun., {\bf 44},
Amer. Math. Soc., Providence, RI.
\bibitem{Birkel-1988}
Birkel, T. (1988). On the convergence rate in the central limit
theorem for associated processes. {\it Ann. Probab.} {\bf 16},
1685--1698.
\bibitem{Bougerol-Picard-1992-a}
Bougerol, P. and Picard, N. (1992). Strict stationarity of
generalized autoregressive processes. {\it Ann. Probab.} {\bf 20},
1714--1730.
\bibitem{Csorgo-1983-LN}
Cs\"{o}rg\H{o}, M. (1983). {\it Quantile Processes with Statistical
Applications}. CBMS-NSF Regional Conference Series in Applied
Mathematics.
\bibitem{Csorgo-Revesz-1978}
Cs\"{o}rg\H{o}, M. and R\'{e}v\'{e}sz, P. (1978). Strong
approximation of the quantile process. {\it Ann. Statist.} {\bf 6},
882--894.
\bibitem{Csorgo-Szyszkowicz-1998}
Cs\"{o}rg\H{o}, M. and Szyszkowicz, B. (1998). Sequential quantile
and Bahadur-Kiefer processes. {\it Order statistics: theory \&
methods}, 631--688, Handbook of Statist., 16, North-Holland,
Amsterdam.
\bibitem{Deheuvels-Mason-1990}
Deheuvels, P. and Mason, D. M. (1990). Bahadur-Kiefer-type
processes. {\it Ann. Probab.} {\bf 18}, 669--697.
\bibitem{Doukhan-1984}
Doukhan, P. (1984). {\it Mixing: Properties and Examples}. Lecture
notes in Statisitcs. Springer-Verlag.
\bibitem{Engle-1982}
Engle, R.F. (1982). Autoregressive conditional heteroscedasticity
with estimates of the variance of United Kingdom inflation. {\it
Econometrica} {\bf 50}, 987--1007.
\bibitem{Esary-Proschan-Walkup-1967}
Esary, J.D., Proschan, F. and Walkup, D.W. (1967). Association of
random variables, with applications. {\it Ann. Math. Statist.} {\bf
38}, 1466--1474.
\bibitem{Fotopoulos-Ahn-1994}
Fotopoulos, S.B. and Ahn, S.K. (1994). Strong Approximation of the
Qunatile Processes and Its Applications under Strong Mixing
Properties. {\it J. Mult. Anal.} {\bf 51}, 17--45.
\bibitem{Ho-Hsing-1997}
Ho, H.-C. and Hsing, T. (1997). Limit theorems for functionals of
moving averages. {\it Ann. Probab.} {\bf 25}, 1636--1669.
\bibitem{Hsing-Wu-2004}
Hsing, T. and Wu, W.B. (2004) On weighted $U$-statisitcs for
stationary sequences. {\it Ann. Probab.} {\bf 32}, 1600--1631.
\bibitem{Kazakevicius-Leipus}
Kazakevicius, V. and Leipus, R. (2003). A new theorem on the
existence of invariant distribution with applications to ARCH
processes. {\it J. Appl. Probab.} {\bf 40}, 147--162.
\bibitem{Kiefer-1967}
Kiefer, J. (1967). On Bahadur's representation of sample qunatiles.
{\it Ann. Math. Statist.} {\bf 38}, 1323--1342.
\bibitem{Kiefer-1970}
Kiefer, J. (1970). Deviations between the sample qunatile process
and the sample df. {\it nonparametric Techniques in Statistical
Inference}, 349--357, M.L. Puri, ed., Cambridge University Press.
\bibitem{Petrov-1975}
Petrov, V.V. (1975). {\it Sums of Independent Random Variables}.
Springer-Verlag.
\bibitem{Wu-2005}
Wu, W.B. (2005). On the Bahadur representation of sample quantiles
for dependent sequences. {\it Ann. Statist.} {\bf 33}, 1934--1963.
\bibitem{Wu-Mielniczuk-2002}
Wu, W. B. and Mielniczuk, J. (2002). Kernel density estimation for
linear processes. {\it Ann. Statist.} {\bf 30}, 1441--1459.
\bibitem{Yu-1996}
Yu, H. (1996). A note on strong approximation for quantile processes
of strong mixing sequences. {\it Stat. Prob. Letters.} {\bf 30},
1--7.
\bibitem{Yu-1995}
Yu, H. (1996). A strong invariance principle for associated
sequences. {\it Ann. Probab.} {\bf 24}, 2079--2097.
\end{thebibliography}
\end{document}